\documentclass[12pt,a4paper]{amsart}
\usepackage{graphics}
\usepackage{epsfig}
\usepackage{graphicx}

\newcommand{\bdi}{\begin{diagram}}
\newcommand{\edi}{\end{diagram}}

\usepackage[colorlinks]{hyperref}
\usepackage{enumerate, amssymb}
\advance\hoffset-10mm \advance\textwidth40mm

\newtheorem{teo}{Theorem}{\bf}{\it}
\newtheorem{conj}[teo]{Conjecture}{\bf}{\it}
\newtheorem{prop}[teo]{Proposition}{\bf}{\it}
\theoremstyle{plain}

\newcommand{\Spec}{ \operatorname{{\rm Spec}}}

\newcommand{\Aut}{ \operatorname{{\rm Aut}}}
\newcommand{\Mat}{ \operatorname{{\rm Mat}}}

\newcommand{\SL}{ \operatorname{{\rm SL}}}

\newcommand{\Susp}{ \operatorname{{\rm Susp}}}

\newcommand{\C}{\ensuremath{\mathbb{C}}}

\newcommand{\Q}{\ensuremath{\mathbb{Q}}}

\newcommand{\G}{\ensuremath{\mathbb{G}}}

\newcommand{\kk}[1]{\bk^{[#1]}}

\newcommand{\cB}{{\ensuremath{\mathcal{B}}}}

\newcommand{\cG}{{\ensuremath{\mathcal{G}}}}

\newcommand{\p}{\partial}

\def\bals#1\eals{\begin{align*}#1\end{align*}}
\def\bal#1\eal{\begin{align}#1\end{align}}

\def\SAut{\mathop{\rm SAut}}
\def\kk{{\Bbbk}}
\def\AA{{\mathbb A}}

\def\ZZ{{\mathbb Z}}

\def\kk{{\Bbbk}}

\def\PP{{\mathbb P}}

\def\sl{\mathfrak{sl}}
\def\reg{{\mathop{\rm reg}}}

\def\and{\quad\mbox{and}\quad}

\renewcommand{\epsilon}{\varepsilon}
\renewcommand{\phi}{\varphi}

\newcommand{\bnum}{\begin{enumerate}}
\newcommand{\enum}{\end{enumerate}}
\renewcommand{\emptyset}{\varnothing}
\newcommand{\brem}{\begin{rem}}
\newcommand{\brems}{\begin{rems}}
\newcommand{\erem}{\end{rem}}
\newcommand{\erems}{\end{rems}}
\newcommand{\bprob}{\begin{prob}}
\newcommand{\eprob}{\end{prob}}
\newcommand{\bprobs}{\begin{probs}}
\newcommand{\eprobs}{\end{probs}}
\newcommand{\bques}{\begin{ques}}
\newcommand{\eques}{\end{ques}}
\newcommand{\bexa}{\begin{exa}}
\newcommand{\bexas}{\begin{exas}}
\newcommand{\eexa}{\end{exa}}
\newcommand{\eexas}{\end{exas}}
\newcommand{\bdefi}{\begin{defi}}
\newcommand{\edefi}{\end{defi}}
\newcommand{\bdefis}{\begin{defis}}
\newcommand{\edefis}{\end{defis}}
\newcommand{\bcor}{\begin{cor}}
\newcommand{\ecor}{\end{cor}}
\newcommand{\blem}{\begin{lem}}
\newcommand{\elem}{\end{lem}}
\newcommand{\bconv}{\begin{conv}}
\newcommand{\econv}{\end{conv}}
\newcommand{\bconj}{\begin{conj}}
\newcommand{\econj}{\end{conj}}
\newcommand{\bprop}{\begin{prop}}
\newcommand{\eprop}{\end{prop}}
\newcommand{\bthm}{\begin{thm}}
\newcommand{\ethm}{\end{thm}}
\newcommand{\bnota}{\begin{nota}}
\newcommand{\enota}{\end{nota}}
\newcommand{\bsit}{\begin{sit}}
\newcommand{\esit}{\end{sit}}
\newcommand{\be}{\begin{equation}}
\newcommand{\ee}{\end{equation}}
\newcommand{\bproof}{\begin{proof}}
\newcommand{\eproof}{\end{proof}}
\def\ba{\begin{array}}
\def\ea{\end{array}}

\begin{document}
\title[Infinite transitivity on affine varieties]
{Infinite transitivity on affine varieties }

\author{I.\ Arzhantsev, H.\ Flenner, S.\ Kaliman, F.\ Kutzschebauch,
M.\ Zaidenberg}
\address{Department of Algebra, Faculty of Mechanics and Mathematics,
Moscow State University, Leninskie Gory 1, GSP-1, Moscow, 119991,
Russia } \email{arjantse@mccme.ru}
\address{Fakult\"at f\"ur
Mathematik, Ruhr Universit\"at Bochum, Geb.\ NA 2/72,
Universit\"ats\-str.\ 150, 44780 Bochum, Germany}
\email{Hubert.Flenner@rub.de}
\address{Department of Mathematics,
University of Miami, Coral Gables, FL 33124, USA}
\email{kaliman@math.miami.edu}
%\author{Frank Kutzschebauch}
\address{Mathematisches Institut\\
Universit\"at Bern\\
Sidlerstrasse 5, CH-3012 Bern, Switzerland, }
\email{frank.kutzschebauch@math.unibe.ch}
\address{Universit\'e Grenoble I, Institut Fourier, UMR 5582
CNRS-UJF, BP 74, 38402 St. Martin d'H\`eres c\'edex, France}
\email{Mikhail.Zaidenberg@ujf-grenoble.fr}

\begin{abstract}
In this note we survey recent results on automorphisms of
affine algebraic varieties, infinitely transitive group actions
and flexibility.
We present related constructions and examples,
and discuss geometric applications and open problems.
\end{abstract}

\maketitle
\date{today}

\thanks{
%\mbox{\hspace{11pt}}{\it 2010 Mathematics Subject Classification}:
%14R20, 14L30.\\
{\renewcommand{\thefootnote}{} \footnotetext{ 2010
\textit{Mathematics Subject Classification:}
14R20,\,32M17; secondary 14L30.\mbox{\hspace{11pt}}}}}

\vfuzz=2pt

\section*{Introduction}\label{sec:intro}

An action of a group $G$ on a set $A$ is said to be {\it $m$-transitive} if
for every two tuples of pairwise distinct points $(a_1,\ldots,a_m)$
and $(a_1',\ldots,a_m')$ in $A$ there exists $g\in G$ such that
$g\cdot a_i=a_i'$ for all $i=1,\ldots,m$. An action which is $m$-transitive
for all $m\in\ZZ_{>0}$ will be called {\it infinitely transitive}.

Clearly, the group of all bijections of an infinite set $A$ acts
infinitely transitively on~$A$. Infinite transitivity never occurs
if $G$ is a Lie group or algebraic group acting on a variety~$A$.
Indeed, $m$-transitivity implies that the map $G\to A^m$ with
$g\mapsto (g.a_1,\ldots, g.a_m)$ is dominant for an $m$-tuple of
pairwise distinct points $(a_1,\ldots,a_{m})~\in~A^m$. This shows
that $G$ cannot act on $A$ $m$-transitively if $\dim G< m\cdot\dim
A$. According to A.~Borel a much stronger result is valid: a real
Lie group can not act even 3-transitively on a simply connected,
non-compact real manifold (see Theorems~5 and 6 in \cite{Bo}). By
a result of F.~Knop~\cite{Kn}, the most transitive action of an
algebraic group over an algebraically closed field is the
3-transitive action of the group $\text{PSL}_2$ on the projective
line $\PP^1$.

At the same time, the group $\Aut(\AA^n)$ of all algebraic
automorphisms of the affine space $\AA^n$ over an infinite field
acts infinitely
transitively on $\AA^n$ for $n\ge 2$. To obtain this result, it suffices
to use linear  automorphisms
and triangular automorphisms of the form
$$
(x_1,\ldots,x_{n-1},x_n) \, \mapsto \, (x_1,\ldots,x_{n-1},
x_n+P(x_1,\ldots,x_{n-1})),
$$
where $P(x_1,\ldots,x_{n-1})$ is an arbitrary polynomial.
These  automorphisms generate the tame automorphism group T$\Aut(\AA^n)$,
which acts infinitely
transitively on $\AA^n$ for $n\ge 2$.

It was shown in~\cite{AFKKZ}  that for certain (infinite dimensional) groups
of automorphisms
of affine varieties transitivity implies infinite transitivity. We do not try
to present here the results of~\cite{AFKKZ} in full generality, but rather
to concentrate on the most interesting features.

\section{Main results}\label{sec:mr}

Let $X$ be an algebraic variety over a field $\kk$. Unless we explicitly
precise the opposite,
we assume usually that $\kk$ is algebraically closed of characteristic zero.
Consider a regular action $\G_a\times X \to X$ of the
additive group $\G_a=(\kk,+)$ of the ground field on~$X$. The image, say,  $H$
of $\G_a$ in the automorphism group $\Aut(X)$ is a one-parameter
unipotent subgroup of $\Aut(X)$. We let $\SAut(X)$ denote the subgroup of~$\Aut(X)$
generated by all its one-parameter unipotent subgroups. Automorphisms
from the group $\SAut(X)$ will be called {\it special}. In general,
$\SAut(X)$ is a normal subgroup of~$\Aut(X)$.

Denote by $X_{\reg}$ the smooth locus of an algebraic variety $X$.
We say that a point
$x\in X_{\reg}$ is {\it flexible} if the tangent space $T_xX$ is
spanned by the tangent
vectors to the orbits $H\cdot x$ over all one-parameter unipotent
subgroups $H$ in $\Aut(X)$.
The variety $X$ is {\it flexible} if every point $x\in X_{\reg}$ is.
Clearly, $X$
is flexible if one point of $X_{\reg}$ is and the group $\Aut(X)$
acts transitively
on $X_{\reg}$.

The following result conjectured in~\cite{AKZ} is proven
in~\cite[Theorem~0.1]{AFKKZ}.

\begin{teo}\label{main}
Let $X$ be an irreducible affine variety of dimension $\ge 2$.
Then the following
conditions are equivalent.
\begin{enumerate}
\item
The group $\SAut(X)$ act transitively on $X_{\reg}$.
\item
The group $\SAut(X)$ act infinitely transitively on $X_{\reg}$.
\item
The variety $X$ is flexible.
\end{enumerate}
\end{teo}

\section{Examples of flexible varieties}\label{sec:cfv}

We are going to show that the equivalent conditions of Theorem~\ref{main}
are satisfied for wide classes of affine varieties.

\smallskip

{\bf 2.1. Suspensions.}\ Let $X$ be an affine variety. Given a nonconstant
regular function $f\in\kk[X]$, we define a new affine variety
$$
\Susp(X,f)=\{uv-f(x)=0\}\subseteq\AA^2 \times X
$$
called a {\it suspension} over $X$. It is shown in \cite[Theorem~3.2]{AKZ}
that
{\em a suspension over a flexible affine variety is again flexible}.
The case of suspensions
over affine spaces was treated earlier in~\cite{KZ1}. Iterating the
construction of
suspension yields new examples of flexible varieties.

Flexibility and infinite transitivity of the action of $\SAut(X)$ is established
in~\cite[Theorem~3.1]{AKZ} for a suspension $X=\{uv-f(x)=0\}$ over
the affine line $\AA^1$
under the assumption that $f(\kk)=\kk$, where $\kk$ is an arbitrary field
of characteristic zero.
The same  holds
for suspensions over flexible real affine algebraic varieties with
connected smooth loci
(\cite[Theorem~3.3]{AKZ}).
By \cite{KuMa}, infinite transitivity holds on every connected component
of the smooth loci
of suspensions over flexible real affine varieties.

\smallskip

{\bf 2.2. Affine toric varieties.}\ Recall that a normal algebraic variety $X$ is
{\it toric}
if it admits a regular action of an algebraic torus $T$ with an open orbit.
In general,
an affine toric variety does not need to be flexible. For instance, if $X=T$ then
the algebra $\kk[X]$ is generated by invertible functions and hence
the group $\SAut(X)$ is trivial.

We say that an affine toric variety $X$ is {\it nondegenerate} if the
only invertible
regular functions on $X$ are nonzero constants. Equivalently, $X$ is
nondegenerate if
it is not isomorphic to $X'\times (\AA^1\setminus \{0\})$ for some
toric variety $X'$. By \cite[Theorem~2.1]{AKZ} any nondegenerate affine
toric variety
is flexible. Considering affine toric surfaces, one obtains examples of
affine varieties $X$
such that $X_{\reg}$ is not a homogeneous space of an algebraic group,
but the group
$\SAut(X)$ acts on $X_{\reg}$ infinitely transitively, see \cite[Example~2.2]{AKZ}.

\smallskip

{\bf 2.3. Homogeneous spaces.}\  Let us consider (following~\cite{Po2})
the class of
connected linear algebraic groups $G$ generated by one-parameter unipotent
subgroups. A connected linear algebraic group $G$ belongs to
 this class if and only if $G$ does not admit non-trivial characters,
 or, equivalently,
if a maximal reductive subgroup of $G$ is semisimple. If such a group $G$
acts on a variety $X$ then the image of $G$ in $\Aut(X)$ is contained in
$\SAut(X)$.
If $G$ acts on $X_{\reg}$ transitively then $X$ is flexible.

As an example, consider a simple rational $G$-module $V$, where $G$ is
semisimple.
The cone $X$ of highest weight vectors in $V$ consists of two $G$-orbits,
namely,  the open orbit
$X\setminus\{0\}$ and the origin $\{0\}$ (\cite{PV1}). If $X\neq V$
then $G$ acts on $X_{\reg}$
transitively, hence
the group $\SAut(X)$ is infinitely transitive
on $X_{\reg}$. Note that $X$ may be considered as a (normal)
affine cone over the
flag variety $G/P$, where a parabolic subgroup $P$ is the stabilizer
of a point
in the projectivization $\PP(X)$ of the cone $X$ in $\PP(V)$. In these terms
infinite transitivity for $X$ was proven in \cite[Theorem~1.1]{AKZ}.

Any affine homogeneous space $G/H$ of dimension $\ge 2$ satisfies
the equivalent
conditions of Theorem~\ref{main} provided that $G$ does not admit
non-trivial characters,
see~\cite[Proposition~5.4]{AFKKZ}. In particular,
for any semisimple group $G$ and a reductive subgroup $H\subseteq G$ the
homogeneous space $X=G/H$ is flexible and the group $\SAut(X)$
is infinitely transitive on $X$.
This applies as well to $X=G$.

\smallskip

{\bf 2.4. Almost homogeneous varieties.}\
Suppose that a connected semisimple algebraic group  $G$
acts with an open orbit on an irreducible affine variety $X$.
In this case we say that $X$
is {\em almost homogeneous}.
It turns out that under some
additional assumptions this implies flexibility of $X$.

\smallskip

{\it 2.4.1. The smooth case.}\ Assume that an almost homogeneous
affine variety $X$ is smooth.
Using Luna's
\'Etale Slice Theorem we show in \cite[Theorem~5.6]{AFKKZ} that $X$
is homogeneous
under the action of a semidirect product $G\ltimes V$,
where $V$ is a certain finite-dimensional
$G$-module. In particular $X$ is flexible.

\smallskip

{\it 2.4.2. $\SL_2$-embeddings.}\ Let the group $\SL_2=\SL_2(\kk)$
act with an open orbit on a normal affine threefold $X$. All such
$\SL_2$-threefolds were classified in \cite{Po5}. If $X$ is smooth
then it is flexible by the above argument. For a  singular  $X$
the complement of the open $\SL_2$-orbit consists of a
two-dimensional orbit, say, $O$ and a singular fixed point $p\in
X$.

It is shown in \cite{BH} that $X$ can be obtained as the quotient of an affine
hypersurface $x_0^b=x_1x_4-x_2x_3$ under an action of a
one-dimensional diagonalizable group.
Such a hypersurface is a suspension over $\AA^3$. Using this one can join
a point in $O$ with a point in the open $\SL_2$-orbit by
a $\G_a$-orbit on $X$ and thus to gain flexibility of $X$;
see \cite[Theorem~5.7]{AFKKZ}
for details.

\smallskip

{\bf 2.5. Vector bundles.}\
Let $\pi\colon E\to X$ be a  reduced, irreducible linear space over a flexible
variety $X$, which is a vector bundle over $X_{\reg}$.
Assume that there is an action
of the group $\SAut(X)$ on $E$ such that the action of every
one-parameter unipotent
subgroup is algebraic and the morphism $\pi$ is equivariant. It is shown in
\cite[Corollary~4.5]{AFKKZ} that the total space $E$ is a flexible variety.
In particular, the tangent bundle $TX$ and all tensor bundles
$E=(TX)^{\otimes a}\otimes(T^*X)^{\otimes b}$ are flexible.

\smallskip

{\bf 2.6. Affine cones over projective varieties.}\
Let $X$ be the affine cone over a projective variety $Y$  polarized by a
very ample divisor $H$. Then one can characterize
flexibility of $X$ in terms of certain geometric properties
of the pair $(Y,H)$ as follows (see
\cite{Pe}).

An open subset $U\subseteq Y$ is called a {\it cylinder} if $U\cong Z\times\AA^1$,
where $Z$ is a smooth affine variety (see \cite{KPZ}, \cite{KPZ1}).
% with $\Pic Z = 0$.
A cylinder $U$ is called {\it $H$-polar} if $U = Y \setminus \text{Supp}\,D$
for some effective $\Q$-divisor $D$ linearly equivalent to $H$.
% with some $d > 0$.
It is shown in \cite[Theorem~3.9]{KPZ} that any $H$-polar
cylinder $U$ on $Y$ gives rise to
a $\G_a$-action on the affine cone $X$ over $Y$.

A subset $W\subseteq Y$ is called {\it invariant} with respect to
a cylinder
$U \cong Z \times \AA^1$  if $W \cap U = \pi^{-1}(\pi(W))$,
where $\pi\colon  U \to Z$ is the first projection. In other
words, $W$ is invariant if every $\AA^1$-fiber of the cylinder is
either contained in $W$
or does not meet~$W$.
A variety $Y$ is {\it transversally covered} by cylinders $U_i$,
$i=1,\ldots,s$,
if $Y = \bigcup_i U_i$ and there is no proper subset $W\subseteq Y$
invariant with respect
to all the $U_i$.

Theorem~2.5 in \cite{Pe} states that
if for some very ample divisor $H$ on a normal projective variety $Y$
there exists
a transversal covering by $H$-polar cylinders then the corresponding
affine cone $X$ over
$Y$ is flexible.
This criterion allows to establish that any affine cone over a
del Pezzo surface of degree
$\ge$ 5 is flexible. The same is true for certain affine cones over del Pezzo
surfaces of degree 4,
including the pluri-anticanonical ones.
In contrast, the pluri-anticanonical cones over del Pezzo surfaces of degree
1 or 2 do not admit
any non-trivial action of a unipotent algebraic group, neither any effective
action of a
two-dimensional connected  algebraic group (\cite{KPZ1}).
The case of cubic surfaces remains open.

\smallskip

{\bf 2.7. Gizatullin surfaces.} These are normal affine surfaces
which admit a completion by a chain of smooth rational curves. It
follows from Gizatullin's Theorem (\cite[Theorems 2 and 3]{Gi},
see also \cite{Du}) that a normal affine surface $X$ different
from $\AA^1\times (\AA^1\setminus\{0\})$ is Gizatullin if and only
if the special automorphism group $\SAut(X)$ has an open orbit;
then this open orbit necessarily has a finite complement in $X$.
It was conjectured in \cite{Gi} that if the base field $\kk$ has
characteristic zero then the open $\SAut(X)$-orbit coincides with
$X_{\reg}$ i.e. that

\smallskip

{\it every Gizatullin surface is flexible.}

\smallskip

This is definitely not true in a positive characteristic, where
the automorphism group $\Aut(X)$ of a Gizatullin surface $X$ can
have fixed points that are smooth points of $X$ \cite{DG}. We have
seen in 2.1 that Gizatullin's Conjecture is true for the
Gizatullin surfaces given in $\AA^3$ by equations $xy - f(z) = 0$,
since these are suspensions over the affine line. Yet another class
of flexible Gizatullin surfaces consists of the Danilov-Gizatullin
surfaces, see \cite{Giz}. Recently S.\ Kovalenko constructed a
counterexample to the Gizatullin Conjecture over $\C$ ({\em
unpublished}).

We refer the reader to \cite{FKZ} and the references therein for a
study of one-parameter groups acting on Gizatullin surfaces.

\smallskip

\section{Technical tools}\label{sec:tt}

We do not try to expose the proof of Theorem~\ref{main} in detail.
In this section we just present a couple of technical tools which
play a crucial role in the proof. The first one is the well known
correspondence between regular $\G_a$-actions on an affine variety
$X$ and locally nilpotent derivations of the algebra $A=\kk[X]$ of
regular functions on $X$.

\smallskip

{\bf 3.1. Locally nilpotent derivations and their replicas.}\
A derivation $\partial$ of an algebra $A$ is called {\it locally
nilpotent} if for any $a\in A$ there exists $m\in\ZZ_{>0}$ such that
$\partial^m(a)=0$. If the group $\G_a$ acts on $X = \Spec\,A$ then
the associated
derivation $\partial$ of A is locally nilpotent.
It is immediate that for every $f\in\ker\partial$ the derivation
$f\partial$ is again
locally nilpotent.

Conversely, given a locally nilpotent derivation
$\partial\colon A \to A $ and $t\in\kk$, the
map $\exp(t\partial)\colon A\to A$ is an automorphism of A.
Furthermore for $\partial \ne 0$,
$H=\exp(t\partial)$ is a one-parameter unipotent subgroup of $\Aut(A)$.
%isomorphic to $\G_a$.
Via the isomorphism $\Aut(A)\cong\Aut(X)$ given by $g\to(g^{-1})^*$
this yields a one-parameter unipotent subgroup of $\Aut(X)$,
which we denote by the same letter $H$. We refer to~\cite{Fre}
for more details on locally nilpotent derivations.

The algebra of invariants $\kk[X]^H=\ker\partial$ has transcendence degree
$\dim X-1$ over $\kk$.
Given an invariant $f\in\kk[X]^H$ the one-parameter unipotent
subgroup $H_f=\exp(\kk f\partial)$, called a {\it replica} of $H$,
plays an important role in the sequel.
The $H_f$-action has the same general orbits as the $H$-action. However, the zero
locus of $f$ remains pointwise fixed under the $H_f$ -action.
So given a finite set of points chosen on distinct general $H$-orbits
one can find a replica $H_f$ of $H$ that moves all the points but a given one.
If we have at our disposal enough $\G_a$-actions in transversal
directions on $X$ then
by changing the velocity along the corresponding orbits as above, we can
move the given ordered finite set in $X_\reg$ into a prescribed position.
This gives the infinite transitivity of the $\SAut(X)$-action on $X_\reg$.

Let us illustrate the notions of a replica
and of a special automorphism in the case of
an affine space $\AA^n$ over $\kk$.
The group  $\SAut (\AA^n)$
contains the one-parameter
unipotent subgroup of translations in any given direction. The
infinitesimal generator of such a subgroup is a directional
partial derivative. Such a derivative defines a locally nilpotent
derivation of the polynomial ring in $n$ variables, whose
phase flow is the group of translations in this
direction. Its replicas are the one-parameter groups of shears in
the same direction.

As another example, consider the locally nilpotent derivation
$\p=X\frac{\p}{\p Y}+Y\frac{\p}{\p Z}$ of the polynomial ring
$\kk[X,Y,Z]$ and an invariant function $f=Y^2-2XZ\in\ker\p$. The
corresponding replica $H_f$ contains in particular the famous
Nagata automorphism $H_f(1)=\exp(f\cdot\p)\in  \SAut (\AA^3)$,
which is known to be wild; see~\cite{US}.

Notice that any automorphism $\alpha\in\SAut (\AA^n)$ preserves
the usual volume form on $\AA^n$.
Hence $\SAut(\AA^n)\subseteq G_n$, where $G_n$
denotes the subgroup of $\Aut (\AA^n)$ consisting of all
automorphisms with Jacobian determinant $1$. The problem whether
the subgroup $\SAut(\AA^n)$ coincides with $G_n$ is widely open.
Recall that this
is the case in dimension 2 due to the Jung-van der Kulk Theorem.

\smallskip

{\bf Algebraically generated groups.}\
Our second tool is a technique to work with infinite dimensional groups.
We say that a subgroup $H$ of the automorphism group
$\Aut(X)$ is {\it algebraic} if $H$ has a structure of an algebraic group
such that the natural action $H\times X \to X$ is a morphism. A subgroup
$G$ of $\Aut(X)$ is called {\it algebraically generated} if it is
generated as an abstract group by a family $\mathcal{G}$ of connected
algebraic subgroups of $\Aut(X)$. Similar notions were studied in the literature
earlier, see e.g. \cite{Ra}, \cite{Sh}, and more recently \cite{Po2}.

In \cite{AFKKZ} we extend some standard facts on finite-dimensional
algebraic transformation
groups to the case of algebraically generated groups. It is not difficult
to show that
for any point $x\in X$ the orbit $G\cdot x$ is locally closed.
What is more surprising,
one can find (not necessarily distinct) subgroups
$H_1,\ldots,H_s\in\mathcal{G}$ such that
$$
G. x =(H_1\cdot H_2\cdot \ldots\cdot H_s). x
$$
for any $x\in X$, see \cite[Proposition~1.5]{AFKKZ}.

In our setting we obtain the following version of Kleiman's Transversality Theorem
\cite[Theorem~1.15]{AFKKZ}.

\begin{teo} \label{5.40}
Let a subgroup $G\subseteq \Aut(X)$ be
algebraically generated by a system $\mathcal{G}$ of connected algebraic
subgroups closed under conjugation in $G$. Suppose that $G$ acts
with an open orbit $O\subseteq X$.
Then there exist subgroups $H_1,\ldots, H_s\in \cG$ such that for
any locally closed reduced subschemes $Y$ and $Z$ in $O$ one can
find a Zariski dense open subset $U=U(Y,Z)\subseteq H_1\times
\ldots \times H_s$ such that every element $(h_1,\ldots, h_s)\in
U$ satisfies the following.
\bnum[(a)]
\item  The translate $(h_1\cdot\ldots\cdot h_s).Z_\reg$
meets $Y_\reg$
transversally.
\item $\dim (Y\cap (h_1\cdot\ldots\cdot h_s).Z)\le
\dim Y+\dim Z-\dim X$.\\
In particular $Y\cap (h_1\cdot\ldots\cdot h_s).Z=\emptyset$
if $\dim Y+\dim Z<\dim X$.
\enum
\end{teo}

The next generalization concerns the Rosenlicht Theorem on rational invariants.
It turns out that for any algebraically generated subgroup $G\subseteq\Aut(X)$
there exists
a finite collection of rational $G$-invariants on $X$ which separate
$G$-orbits in general position \cite[Theorem~1.13]{AFKKZ}. In particular,
the codimension of a general $G$-orbit in $X$ equals the transcendence degree
of the field $\kk(X)^G$ of rational $G$-invariants over $\kk$.
The latter result has a useful corollary.

\smallskip

{\bf The Makar-Limanov invariant.} \ Recall~\cite{Fre} that the
{\it Makar-Limanov invariant} $\text{ML}(X)$
of an affine algebraic variety $X$ is the
intersection of the kernels of all locally nilpotent derivations on $\kk[X]$.
In other words, $\text{ML}(X)$ is the subalgebra of all $\SAut(X)$-invariants
of the algebra $\kk[X]$.
Similarly~\cite{Lie1} the {\it field Makar-Limanov invariant} $\text{FML}(X)$
is defined as the intersection of the kernels of extensions
of all locally nilpotent derivations on $\kk[X]$ to the field of fractions $\kk(X)$.
This is a subfield of $\kk(X)$ which consists of all
rational SAut(X)-invariants. If it is trivial i.e., if
$\text{FML}(X)=\kk$, then so is $\text{ML}(X)$, while the converse is not true in general.
Triviality of  $\text{FML}(X)$ is equivalent to the existence of a
flexible point in $X_\reg$,
and to the existence of an open $\SAut(X)$-orbit in $X$.

The question arises how these invariants are connected with rationality
properties of the variety
$X$. There are examples of non-unirational affine threefolds X with
$\text{ML}(X)=\kk$ birationally
equivalent to $C\times\AA^2$, where $C$ is a curve of genus $g\ge 1$,
see~\cite[Example~4.2]{Lie2}. For such a
threefold $X$ the general $\SAut(X)$-orbits have dimension two,
the field Makar-Limanov invariant $\text{FML}(X)$ is non-trivial,
and there is no
flexible point in $X$.

The next proposition confirms, in particular, Conjecture~5.3 in~\cite{Lie1}
(cf.\ also~\cite{BKK} and~\cite{Po2}).

\begin{prop} \cite[Proposition~5.1]{AFKKZ}
Let $X$ be an irreducible affine variety. If the field Makar-Limanov invariant
$\text{FML}(X)$ is trivial then $X$ is unirational.
\end{prop}

Indeed, the condition
$\text{FML}(X)=\kk$ implies that the group $\SAut(X)$ acts on $X$ with an
open orbit $O$. Thus  there are $\G_a$-subgroups $H_1,\ldots,H_s$ in $\SAut(X)$
and a point $x\in X$ such that the image of the map
$$
H_1\times\ldots\times H_s \to X, \quad (h_1,\ldots,h_s)
\mapsto (h_1\ldots h_s). x
$$
coincides with $O$. Since $H_1\times\ldots\times H_s$
is isomorphic (as a variety) to the affine space $\AA^s$,
this yields unirationality of $X$.
Moreover, any two points in
$O$ are contained in the image of a morphism $\AA^1\to O$.
In particular, $O$
is $\AA^1$-connected in the sense of \cite[6.2]{KK2}.

In general, flexibility implies neither rationality nor stable rationality.
Indeed,
there exists a finite subgroup $F\subset\SL_n$, where $n\ge 4$,
such that the smooth
unirational affine variety $X=\SL_n/F$ is not stably rational, see
\cite[Example~1.22]{Po2}. However, by 2.3 the variety $X$
is flexible and the group
$\SAut(X)$ acts infinitely transitively on $X$.

We expect further development of the invariant theory
for algebraically generated groups.

\section{Geometric consequences}\label{sec:gc}

Let us start with several results related to Theorem~\ref{main}.

\smallskip

{\bf 4.1. Collective transitivity.}\
By a {\it collective infinite transitivity} we
mean a possibility to move simultaneously (that is,
by the same automorphism)
an arbitrary finite set of points along their orbits into
a given position.
We illustrate our general results in this direction on a concrete
example from
linear algebra, cf.~\cite{Re}.

Let $X=\Mat(n,m)$ be the space of all $n\times m$ matrices over
$\kk$.
The subset  $X_r\subseteq X$ of matrices of rank $r$ is well known
to have dimension $mn-(m-r)(n-r)$. In the following we always
assume that this dimension is $\ge 2$.
The product $\SL_n\times\SL_m$ acts on $X$ via the
left-right multiplication preserving the strata $X_r$. For every
$k\neq l$ we let $E_{kl}\in\sl_n$ and $E^{kl}\in\sl_m$ denote the
nilpotent matrices with $x_{kl}=1$ and the other entries equal
zero. Let further $H_{kl}=I_n+\kk E_{kl}\subseteq \SL_n$ and
$H^{kl}=I_m+\kk E^{kl}\subseteq \SL_m$ be the corresponding
one-parameter unipotent subgroups acting on the stratification
$X=\bigcup_r X_r$, and let $\delta_{kl}$ and $\delta^{kl}$,
respectively, be the corresponding locally nilpotent vector
fields on $X$ tangent to the strata.

We call {\em elementary} the one-parameter unipotent subgroups
$H_{kl}$, $H^{kl}$, and all their replicas. In the following
theorem we establish the collective infinite transitivity on the
above stratification of the subgroup $G$ of $\SAut(X)$ generated
by the two sides elementary subgroups.

By a well known theorem of linear algebra,
the subgroup $\SL_n\times\SL_m\subseteq G$ acts transitively
on each stratum $X_r$ (and so these strata are $G$-orbits)
except for the open stratum $X_n$
in the case where $m=n$. In the latter case
the $G$-orbits contained in $X_n$
are the level sets of the  determinant.

\begin{teo} \cite[Theorem~3.3]{AFKKZ}
Given two finite ordered collections $\cB$
and $\cB'$ of distinct matrices in $\Mat(n,m)$ of the same
cardinality, with the same sequence of ranks,
and in the case where $m=n$ with the same sequence of
determinants, we can simultaneously transform  $\cB$ into  $\cB'$
by means of an element $g\in G$, where
$G\subseteq\SAut(\Mat(n,m))$ is the subgroup generated by all
elementary one-parameter unipotent subgroups.
\end{teo}

See \cite[Section~3.3]{AFKKZ} for similar results on symmetric and
skew-symmetric matrices.

\smallskip

{\bf 4.2. $\AA^1$-richness.}\
Let $X$ be a flexible affine variety of dimension $\ge 2$,
and let $p_1,\ldots,p_k\in X_{\reg}$ be a $k$-tuple.
Fix a $\G_a$-orbit $C$ on $X$ and some $k$-tuple  of distinct points
$q_1,\ldots,q_k\in C$.
Due to infinite transitivity there is an element $g\in\SAut(X)$ such that
$g\cdot q_1=p_1,\ldots,g\cdot q_k=p_k$. So the translate $g\cdot C$
of $C$ is a $\G_a$-orbit
on $X$ passing through $p_1,\ldots,p_k$. This elementary observation can be
strengthened  in the following way.

An affine variety $X$ is called {\em
$\AA^1$-rich} if for every finite subset $Z$ and every algebraic
subset $Y$ of codimension $\ge 2$ there is a
curve in $X$ isomorphic to the affine line $\AA^1$, which is
disjoint with $Y$
and passes through every point of $Z$
(\cite{KZ2}).

The following result is immediate from the Transversality
Theorem \ref{5.40}.

\begin{teo} \cite[Corollary~4.18]{AFKKZ}
Let $X$ be an affine variety. Suppose that the group $\SAut(X)$
acts with an open orbit $O\subseteq X$. Then for any finite subset
$Z\subseteq O$
and for any closed subset $Y\subseteq X$ of codimension $\ge 2$ with
$Z\cap Y=\emptyset$ there is an orbit $C\cong\AA^1$ of a
$\G_a$-action on $X$ which does not meet $Y$ and passes through
each point of $Z$.
\end{teo}

In the special case where $X=\AA_\C^n$ this
also follows from the Gromov-Winkelmann Theorem~\cite{Wi}
which says that
the group $\Aut (\AA^n\backslash Y)$ acts transitively on
$\AA^n\setminus Y$,
combined with the equivalence of transitivity and infinite
transitivity
of Theorem \ref{main},
which is valid in this setting as well.
More generally, we also show that $C$ as in the theorem
can be chosen to have  prescribed jets at the points of $Z$.

\smallskip

{\bf 4.3. Prescribed jets of automorphisms.}\
Our results on infinite transitivity may be strengthened
in the following way,
see~\cite[Theorem~4.14 and Remark~4.16]{AFKKZ}.

\begin{teo}
Let $X$ be a flexible affine variety of dimension
$n\ge 2$ equipped with an algebraic  volume form\footnote{By this
we mean a nowhere vanishing $n$-form defined on $X_\reg$.}
$\omega$. Then for any $m\ge 0$ and for any finite subset
$Z\subseteq X_\reg$ there exists an automorphism $g \in \SAut(X)$ with
prescribed $m$-jets at the points $p\in Z$, provided these jets
preserve $\omega$ and inject $Z$ into $X_\reg$.
The same holds without the requirement
that there is a global volume form on $X_\reg$
provided that for every $p\in Z$ the
corresponding jet fixes the
point p and its linear part belongs to the group  $\SL(T_pX)$.
\end{teo}

{\bf 4.4. The Oka-Grauert-Gromov Principle for flexible varieties.}
Let us provide an important application of flexibility
in analytic geometry;
see \cite[Theorem 6.2 and Proposition 6.3]{AFKKZ}.
We address the reader to \cite[\S 6]{AFKKZ}
for more details and a survey.

\begin{teo}
Let $\pi : X \to B$ be a surjective submersion of smooth
irreducible affine algebraic
varieties over $\C$ such that for some algebraically generated subgroup
$G \subseteq \Aut(X)$
the orbits of $G$ coincide with the fibers of $\pi$.
Then the Oka-Grauert-Gromov principle
holds for $\pi : X \to B$.  That is, any continuous section of $\pi$
is homotopic to a holomorphic one, and
any two holomorphic sections of $\pi$ that are homotopic
via continuous sections are
also homotopic via holomorphic ones.\end{teo}

\section{Open problems}\label{sec:op}

Let us finish this note with several open problems on flexible varieties.
The examples from subsection~2.4 motivate the following problem.

\smallskip

{\it Characterize flexible varieties among the normal almost homogeneous affine varieties. }

\smallskip

 By the result described in 2.4.1, a {\em smooth} almost homogeneous variety is flexible.
In fact, in all examples that we know an almost homogeneous normal variety is flexible.
For instance,
one might hope for positive results
in the class of spherical varieties. By definition, a $G$-variety $X$ is {\em spherical}
if a Borel subgroup
$B$ of $G$ acts on $X$ with an open orbit.  An important particular case
is the variety $X=\Spec\,\kk[G/U]$, where $U$ is a maximal unipotent subgroup
of a semisimple group $G$.

\smallskip

To formulate the next problem we need to introduce some more notation.
Let $Y$ be a closed subvariety of an affine variety $X$. Denote by
$\SAut(X)_Y$ the subgroup generated by all one-parameter unipotent subgroups
$\exp(\kk\partial)$,
where the locally nilpotent vector field $\partial$ vanishes on $Y$.

\smallskip

{\it Assume that the group $\SAut(X)$ acts on $X$ with an open orbit $O$, and let
$Y\subseteq O$ be a closed subvariety of codimension $\ge 2$. Is it true that
the group $\SAut(X)_Y$ acts on $O\setminus Y$ transitively? In particular, is
$X\setminus Y$ flexible if so is $X$?}

\smallskip

Some positive results on this problem can be found in \cite[Proposition~4.19]{AFKKZ}.

\smallskip

Our last problem concerns exotic structures on the affine spaces.

\smallskip

{\it Does there exist a flexible exotic algebraic
structure on an affine space that is, a flexible smooth affine algebraic
variety over $\C$ diffeomorphic  but not isomorphic to an affine
space $\AA^n_\C$?}

\smallskip

Notice that for all the exotic structures on $\AA^n_\C$ constructed so far,
the Makar-Limanov invariant is non-trivial, whereas for a flexible
such structure even the field
Makar-Limanov invariant must be trivial (cf. however \cite{Du1}).

\end{document}